%Paper: q-alg/9504006
%From: arnaudon@lapphp8.in2p3.fr
%Date: Wed, 12 Apr 1995 10:55:09 +0200
%Date (revised): Wed, 21 Jun 1995 05:32:11 -0400

%23456789-123456789-223456789-323456789-423456789-523456789-623456789-723456789
%---:----1----:----2----:----3----:----4----:----5----:----6----:----7----:6789
\input harvmac.tex
\input epsf
%%%%%%%%%%%%%%%%%%%%%%%%%%%%%%%%%%%%%%%%%%%%%%%%%%%%%%%%%%%%%%%%%%%%%%%%%%%%%%%
%
%    Title: Representations of U_q(sl(N)) at Roots of Unity
%
%    Authors: B. Abdesselam, D. Arnaudon, A. Chakrabarti
%
%    Comments:  Plain TeX, Macros  harvmac.tex  and epsf needed
%               3 figures in a uuencoded tar separate file
%
%               Save the figure file as uqslnfig.uu
%               Type    uudecode uqslnfig.uu
%                       tar xvf uqslnfig.tar
%               Then you can TeX the main file
%
%%%%%%%%%%%%%%%%%%%%%%%%%%%%%%%%%%%%%%%%%%%%%%%%%%%%%%%%%%%%%%%%%%%%%%%%%%%%%%%
%
%    References:
%
%
%\def\bibitem#1#2{\expandafter\lref\csname r#1\endcsname{#2}}
\def\bibitem#1#2\par{\expandafter\lref\csname r#1\endcsname{#2}}

\bibitem {AB}  D. Arnaudon and M. Bauer,
  {\sl Polynomial Relations in the Centre of
    ${\cal U}_q(sl(N))$, }
  hep-th/9310030,
  Lett. Math. Phys. {\bf 30} (1994) 251.

\bibitem {ACi} D. Arnaudon and A. Chakrabarti,
  {\sl Periodic and partially periodic representations of $SU(N)_{q}$,}
  Commun. Math. Phys. {\bf 139} (1991) 461.

\bibitem {ACii} D. Arnaudon and A. Chakrabarti,
  {\sl Flat periodic representations of
    $\ \cU_{q}(\cG)$,} Commun. Math. Phys. {\bf 139} (1991) 605.

\bibitem {Asutroisq} D. Arnaudon,
  {\sl Periodic and flat irreducible representations of
    $SU(3)_{q}$,}
  Commun. Math. Phys. {\bf 134} (1990) 523.

\bibitem {Bab} O. Babelon,
  {\sl Representations of the Yang--Baxter algebrae associated to Toda
    field theory,} Nucl. Phys. {\bf B230} [FS10] (1984) 241.

\bibitem {BHV} \v{C}. Burd\'{\i}k, M. Havl\'{\i}\v{c}ek and
  A. Van\v{c}ura,
  {\sl Irreducible highest weight representations of quantum groups
    $\cU_q(gl(n,C))$},
  Commun. Math. Phys. {\bf 148} (1992) 417.

\bibitem {BKMS} V. V. Bazhanov, R. M. Kashaev, V. V. Mangazeev and Yu.
  G. Stroganov, {\sl $({\bf Z_N}\times)^{n-1}$ Generalization of the
    Chiral Potts Model},
  Commun. Math. Phys. {\bf 138} (1991) 393.

\bibitem {BKW} \v{C}. Burd\'{\i}k, R. C. King and T. A. Welsh,
  {\sl The explicit construction of irreducible representations of the
    quantum algebras $\cU_q(sl(n))$,}
  Proceedings of the 3rd Wigner Symposium, Oxford, September 1993,
  Classical and quantum systems, Eds. L. L. Boyle and A. I. Solomon,
  World Sci., Singapore (1993).

\bibitem {ChariPr} V. Chari and A. Pressley, {\sl Fundamental
    representations of quantum groups at roots of 1} .

\bibitem {DJMMc} E. Date, M. Jimbo, K. Miki and  T. Miwa,
  {\sl Generalized chiral Potts models and minimal cyclic
    representations of
    $\ \cU_q \left(\widehat{gl}(n,\CC )\right) $,}
  Commun. Math. Phys. {\bf 137}  (1991) 133.

\bibitem {DJMMb} E. Date, M. Jimbo, K. Miki and  T. Miwa,
  {\sl Cyclic
    representations of $\ \cU _{q}(sl(n+1,\CC ))$ at $q^{N}=1$,}
                                %Preprint RIMS-703
  (1990), Publ. RIMS., Kyoto Univ. {\bf 27} (1991) 347.

\bibitem {DK}  C. De Concini and V.G. Kac,
  {\sl Representations of quantum groups at roots of 1,}
  Progress in Math. {\bf 92} (1990) 471 (Birkh\"auser).

\bibitem {DKP}  C. De Concini, V.G. Kac and C. Procesi,
  {\sl Quantum coadjoint action,}
                                %Preprint Pisa (1991)
  J. A.M.S. {\bf 5} (1992) 151,
  and
  {\sl Some remarkable degenerations of quantum groups,}
  Commun. Math. Phys. {\bf 157} (1993) 405.

\bibitem {Dob} %\bibitem {DobStAndrews} \bibitem {DobSchlossHofen}
  V.K. Dobrev, in Proc. {\sl Int. Group Theory
    Conference,} St Andrews, 1989, Vol. 1, Campbell and Robertson (eds.),
London
  Math. Soc. Lect. Notes Series 159, Cambridge University Press, 1991
  and {\sl Representations of
    quantum groups,} Proc. ``Symmetries in Science
  V: Algebraic Structures, their Representations, Realizations and
  Physical Applications'', Schloss Hofen, Austria,
  (1990), Eds. B.  Gruber, L.C. Biedenharn
  and H.-D. Doebner (Plenum Press, NY, 1991) pp. 93-135.

\bibitem {Dri} V.G. Drinfeld, {\sl Quantum Groups,} Proc.
  Int. Congress of Mathematicians, Berkeley, California, Vol.
  {\bf 1}, Academic Press, New York (1986), 798.

\bibitem {FGP} P. Furlan, A. C. Ganchev and V. B. Petkova,
  {\sl Quantum groups and fusion rules multiplicities,}
  Nucl. Phys. {\bf B343} (1990) 205.

\bibitem {FRT} L. D. Faddeev, N.Yu. Reshetikhin and L.A.
  Takhtajan, {\sl
    Quantization of Lie groups and Lie algebras}, Leningrad Math.
  J.  {\bf 1} (1990) 193.

\bibitem {GasperRahman} G. Gasper and M. Rahman, {\sl Basic
    hypergeometric series,} Cambridge University Press (1990).

\bibitem {Jim} M. Jimbo, {\sl $q$-difference analogue of
    $\ \cU (\cG )$ and the Yang Baxter equation,} Lett. Math. Phys. {\bf 10}
  (1985)  63.

\bibitem{JimGZ} M. Jimbo, Lecture notes in Physics 246, Springer
  (1986), 334.

\bibitem {KerDipl} T. Kerler,
  {\sl Darstellungen der Quantengruppen und Anwendungen,}
  Diplomarbeit, ETH--Zurich, August 1989.

\bibitem {Lus} G. Lusztig,
  {\sl Quantum groups at roots of $1$,} Geom. Ded. {\bf 35} (1990) 89.

\bibitem {LusJAMS} G. Lusztig,
  {\sl Finite dimensional Hopf algebras arising from quantized universal
    enveloping algebras,} Journ. A.M.S. {\bf 3} (1990) 257.

\bibitem {NAKyStoi} Nguyen Anh Ky and N.I. Stoilova, {\sl
    Finite-dimensional representations of the quantum superalgebra
    $\cU_q[gl(2|2)]$ II: nontypical representations at generic $q$,}
  preprint IC/94/352.

\bibitem {Pal} T.D. Palev,
    Funct. Anal. Appl. vol. 21, No.3 p.85 (1987) (English transl),
    Journ. Math. Phys. vol. 29, 2589 (1988),
    Journ. Math. Phys. vol. 30, 1433 (1989).

\bibitem {PalTol} T.D. Palev and V.N. Tolstoy, {\sl Finite dimensional
irreducible representations of the quantum superalgebra
$\cU_q[gl(n|1)]$,} Commun. Math. Phys. {\bf 141} (1991) 549.

\bibitem {PSVdJ} T.D. Palev, N.I. Stoilova and J. Van der Jeugt, {\sl
    Finite-dimensional representations of the quantum superalgebra
    $U_q[gl(n|m)]$ and related $q$-identities,} Commun. Math. Phys. {\bf
    166} (1994) 367.

\bibitem {RA}  P. Roche and  D. Arnaudon,
  {\sl Irreducible representations of the quantum
    analogue of $SU(2)$,}
  Lett. Math. Phys. {\bf 17} (1989) 295.

\bibitem {RosA} M. Rosso,
  {\sl Finite dimensional representations of the quantum analogue of the
    enveloping algebra of a complex simple Lie algebra,}
  Commun. Math. Phys. {\bf 117} (1988) 581.

\bibitem {RosHC} M. Rosso,
  {\sl Analogues de la forme de Killing et du th\'eor\`eme
    d'Harish--Chandra pour les groupes quantiques,}
  Ann. Scient. \'Ec. Norm. Sup., $4^e$ s\'erie, t.23, (1990), 445.

\bibitem {Schnizer} W. A. Schnizer, {\sl Roots of unity:
    representations of quantum groups,}
  Commun. Math. Phys. {\bf 163} (1994) 293.

\bibitem {Skl} E. K. Sklyanin, {\sl Some algebraic
    structures connected with the
    Yang--Baxter equation. Representations of quantum algebras,}
  Funct. Anal. Appl. {\bf 17} (1983) 273.

\bibitem {Ueno} K. Ueno, T. Takebayashi and Y. Shibukawa.
  Lett. Math. Phys. {\bf 18}, 215 (1989).

%%%%%%%%%%%%%%%%%%%%%%%%%%%%%%%%%%%%%%%%%%%%%%%%%%%%%%%%%%%%%%%%%%%%%%%%%%%%%%%
%
% Quelques Macro supplementaires ...
%
%
%  Lettres grasses majuscules
%

%
%  Lettres calligrafiees
%
                    \def\cC{{\cal C}}
                    \def\cF{{\cal F}}
\def\cG{{\cal G}}

\def\cP{{\cal P}}                    
                    \def\cU{{\cal U}}

%
%   Capital roman double letters
%

\def\CC{\rlap {\raise 0.4ex \hbox{$\scriptscriptstyle |$}}
\hskip -0.1em C}
\def\FF{\hbox to 8.33887pt{\rm I\hskip-1.8pt F}}
\def\NN{\hbox to 9.3111pt{\rm I\hskip-1.8pt N}}
\def\PP{\hbox to 8.61664pt{\rm I\hskip-1.8pt P}}
\def\QQ{\rlap {\raise 0.4ex \hbox{$\scriptscriptstyle |$}}
{\hskip -0.1em Q}}
\def\RR{\hbox to 9.1722pt{\rm I\hskip-1.8pt R}}
\def\ZZ{\hbox to 8.2222pt{\rm Z\hskip-4pt \rm Z}} %\def\ZZ{Z\!\!\!Z}
%
%
% Definitions locales
%
\def\uqn{\cU_q(sl(N))}
\def \d2dots{\mathinner{\mkern1mu\raise1pt\vbox{\kern7pt\hbox{.}}\mkern2mu
\raise4pt\hbox{.}\mkern2mu\raise7pt\hbox{.}\mkern1mu}}
\def\fp{\hbox{f.p.}}
\font\titlefonti=cmssbx10 scaled \magstep3
\font\titlefontii=cmsy10 scaled \magstep3
\font\titlefontiii=cmmib10 scaled \magstep3

\def\enslapp{ENSLAPP}
%
%%%%%%%%%%%%%%%%%%%%%%%%%%%%%%%%%%%%%%%%%%%%%%%%%%%%%%%%%%%%%%%%%%%%%%%%%%%%%%
%%%
%%%%%%%%%%%%%%%%%%%%%%%%%%%%%%%%%%%%%%%%
%%%%%%%% LOGO ENSLAPP - DEBUT  %%%%%%%%%
%%%%%%%%%%%%%%%%%%%%%%%%%%%%%%%%%%%%%%%%
\hbox to \hsize{\vbox{\hsize 8 em
{\bf \centerline{Groupe d'Annecy}
\ \par
\centerline{Laboratoire}
\centerline{d'Annecy-le-Vieux de}
\centerline{Physique des Particules}}}
\hfill
%\vbox{\hsize 20 em \hbox{\raise -3 ex\hbox{\epsfbox{
%/lapphp8/keklapp/ragoucy/paper/enslapp.ps}}}}
\hfill
\vbox{\hsize 7 em
{\bf \centerline{Groupe de Lyon}
\ \par
\centerline{Ecole Normale}
\centerline{Sup\'erieure de Lyon}}}}
\hrule height.42mm
\vfill
%%%%%%%%%%%%%%%%%%%%%%%%%%%%%%%%%%%%%%%%
%%%%% LOGO ENSLAPP  - FIN %%%%%%%%%%%%%%
%%%%%%%%%%%%%%%%%%%%%%%%%%%%%%%%%%%%%%%%
%%%
\Title{}
{\titlefonti
{\vbox{\centerline{Representations of
{\titlefontii U}$_q${\titlefontiii $($sl$($N$))$}
%\uqn
at Roots of Unity }}}
}
\centerline{B. Abdesselam
$^{a,}$\footnote{$^{1}$}{abdess@orphee.polytechnique.fr},
 D. Arnaudon $^{b,}$\footnote{$^{2}$}{arnaudon@lapp.in2p3.fr},
A. Chakrabarti$^{a,}$\footnote{$^{3}$}{chakra@orphee.polytechnique.fr}}
\vskip 1truecm
\centerline{\it $^{a}$
Centre de Physique Th{\'e}orique,
Ecole Polytechnique
\footnote{$^4$}
{\it Laboratoire Propre du CNRS UPR A.0014}}
\centerline{\it 91128 Palaiseau Cedex, France.}
\smallskip
\centerline{\it $^{b}$
ENSLAPP \footnote{$^{5}$}
{\it URA 14-36 du CNRS, associ\'ee \`a l'E.N.S. de Lyon et \`a
l'Universit\'e de Savoie.},
Chemin de Bellevue BP 110,}
\centerline{\it 74941 Annecy-le-Vieux Cedex, France.}

\vfill

%%%  ABSTRACT:
The Gelfand--Zetlin basis for representations of $\uqn$ is improved to
fit better the case when $q$ is a root of unity.
The usual $q$-deformed representations, as well
as the nilpotent, periodic (cyclic), semi-periodic (semi-cyclic) and
also some atypical
representations are now described with the same formalism.

\vfill

\rightline{                                                     q-alg/9504006}
\rightline{                                                       RR343.01.95}
\rightline{                                                 \enslapp-A-506/95}
\rightline{March 95}

\baselineskip=15.1pt

\Date{}
%\draft

\newsec{Introduction}
We are interested in
quantum Lie algebras \refs{\rDri, \rJim, \rFRT} and their
finite dimensional irreducible representations.
At generic deformation parameter $q$, the classification of
irreducible representations is in correspondence with the classical
case \rRosA.
When $q$ is a $m$-th
root of unity, there are two options:

One can consider
first the restricted quantum Lie algebra, where the raising and
lowering generators are nilpotent,
i.e. $e_\alpha^{m}=f_\alpha^{m}=0$ and where the Cartan generators
$h_i$ are such that
$k_i^{m}=\left(q^{h_i}\right)^{m}=1$. Representations of those were
studied by Lusztig \rLus.
A classification of irreducible representations of $\cU_q(sl(3))$
was done in \rDob.

The other option is to fix no relation for the $m$-th powers of these
generators, which are actually, for an odd value of $m$,
in the centre of the quantum algebra
\rDK. Then the irreducible representations may admit a periodic action
for the raising and lowering generators. An important work has already
been done towards the classification of these representations
\refs{\rDK,\rDKP}.  It seems to us
that apart from the $\cU_q(sl(2))$ \rRA\
case there is still however no complete classification.

On the other hand, there are already explicit expressions for
representations of $\uqn$ at roots of unity,
the case we will consider from now on.

In \refs{\rSkl \rBab \rAsutroisq \rDJMMb \rDJMMc \rACi \rBKMS \rChariPr
{--} \rSchnizer},
explicit expressions
for representations with periodic (or cyclic) actions of the
generators are given.

In \refs{\rJimGZ \rBHV \rBKW {--} \rUeno}, explicit expressions for usual
representations of $\uqn$ are written, which lead when the deformation
parameter $q$ goes to a root of unity to either irreducible or
reducible (sometimes not totally reducible) representations, depending
on their highest weight \refs{\rDob, \rBKW, \rFGP}.

The irreducible sub-factors of representations that become reducible
or indecomposable in the limit where $q$ is a root of unity have
sometimes no classical counterpart \refs{\rDob,\rBKW}
and we call them
atypical by analogy with the case of representations of superalgebras.
Some of them also appear as sub-factors of some degenerations of periodic
representations \rAsutroisq.

In this paper, we present an improvement  of \rACi\
based on the Gelfand--Zetlin construction, which allows us to write
explicitly, with the same formalism, irreducible representations of
$\uqn$ at roots of unity independently of their nature
(i.e. periodic, semi-periodic, usual or some atypical). Their nature is
actually encoded in the generalized
parameters involved in the Gelfand--Zetlin
basis.
All the types of finite dimensional irreducible representations we are
aware of enter in this scheme (however, atypical representations
generally need a special treatment).

In Sect. 2, we present the formalism and some general rules for the
construction of finite dimensional irreducible representations using
the adapted Gelfand--Zetlin pattern.
The main types of representations are presented as examples in Sect 3.
As an application, we finally
give in Sect. 4
a set of relations among the generators of
the centre of $\uqn$ that generalizes the relations derived in
\refs{\rKerDipl, \rAB}.

\newsec{The adapted Gelfand--Zetlin basis}

\subsec{The quantum algebra $\uqn$}

The quantum
algebra $\uqn$ \refs{\rDri, \rJim}\ is defined by the generators $k_i$,
$k_i^{-1}$, $e_i$, $f_i$ ($i=1,...,N-1$) and the relations
\eqn\euqn{
\eqalign{ k_i e_j & = q^{a_{ij}} e_j k_i
          \qquad \qquad
          k_i f_j  = q^{-a_{ij}} f_j k_i \;,\cr
          [e_i,f_j] & = \delta_{ij} {k_i -k_i^{-1} \over q-q^{-1}} \;,\cr
          [e_i,e_j] & = 0 \qquad \hbox{for} \qquad |i-j|>1 \;,\cr
          e_i^2 e_{i\pm 1} & - (q+q^{-1}) e_i e_{i\pm 1} e_i + e_{i\pm
            1} e_i^2 = 0 \;,\cr
          f_i^2 f_{i\pm 1} & - (q+q^{-1}) f_i f_{i\pm 1} f_i + f_{i\pm
            1} f_i^2 = 0 \;,\cr
}}
where $\left( a_{ij} \right)_{i,j=1,...,N-1}$ is the Cartan matrix of
$sl(N)$, i.e. $a_{ii}=2$, $a_{i,i\pm 1}=-1$ and $a_{ij}=0$ for
$|i-j|>1$.

We will not use the (standard) co-algebra structure in the following.

Let us now define the adapted Gelfand--Zetlin basis for the
representations of $\uqn$.

\subsec{Vectors of the Gelfand--Zetlin basis}

The states are

\eqn\estates{ |p\rangle =
\left|
\matrix{&{p}_{1N}&&{p}_{2N}&\cdots&{p}_{N-1,N}&&{p}_{NN}\cr
\cr
&&{p}_{1N-1}&&\cdots&&{p}_{N-1,N-1}&\cr
\cr
&&&\hskip -.8cm \ddots &\cdots&\hskip .8cm \d2dots &&\cr
\cr
&&&{p}_{12}&&{p}_{22}&&\cr
\cr
&&&&{p}_{11}&&&\cr}
\quad\right.
\hbox{\raise -12 ex\hbox{\epsfbox{uqsln3.eps}}}
}
(with respect to \rACi, we use $p_{il} = h_{il}-i$ instead of
$h_{il}$.)

As usual, the indices of the first line are fixed on a given
representation, whereas the others move by steps of $\pm 1$ under the
action of the raising and lowering generators.
The whole set of $p_{il}$'s is defined up to an overall constant (only
differences enter in formulas). One can constrain for example
$\sum p_{iN}$, or $p_{NN}$, to be zero.
The actions of the
Cartan generators are diagonal on this basis.

In the classical case, and also
in the quantum case when $q$ is generic, the indices $p_{il}$
are integer and satisfy the triangular identities
\eqn\etriangiden{p_{i,l+1} \geq p_{il} > p_{i+1,l+1}\;.}
The first line of indices determines in this case the highest weight
of the representation.

In the case we consider in the following, $q$ is a root of unity, and
the indices $p_{il}$ are {\it complex}. Let $m$ be the smallest integer
such that $q^{m}=1$. We will only consider the case of odd $m$ in this
paper. Most of the result can be adapted to the case where $m$ is even
by applying the prescription of \rACi.

Since the indices $p_{il}$ appear in the expressions of the matrix
elements only through the quantities $q^{p_{il}}$,
they can consistently be defined modulo $m$ for most of the
representations, i.e. two states with indices differing by multiples
of $m$ can be identified. This will be our convention unless we
specify it in the text.

We define the ``fractional part'' $\fp(p_{il})$ of
$p_{il}$ by
\eqn\efracpart{\fp(p_{il}) =
p_{il} \quad [\hbox{mod} \quad {1\over 2}]
}

Since the generators move the indices $p_{il}$ by integer steps,
their fractional part is then fixed on a representation,
even if they do not belong to the first line.

The specifications and restrictions on the values of the indices
$p_{il}$ will be given after the action of the generators on the
vectors $|p\rangle$.

\subsec{Action of the generators}

The action of the generators $k_l^{\pm 1}$, $e_l$, $f_l$ is given by
\eqn\erepr{\eqalign{
k_l^{\pm 1} |p\rangle & =
q^{\pm\left( 2\sum\limits_{i=1}^{l}p_{il}
             -\sum\limits_{i=1}^{l+1}p_{i,l+1}
             - \sum\limits_{i=1}^{l-1}p_{i,l-1} - 1 \right)}
|p\rangle \;,\cr\cr
f_l |p\rangle & = \sum_{j=1}^l
c_{jl} {P'_1(j,l;p) P'_2(j,l;p) \over P'_3(j,l;p) }
|p_{jl}-1 \rangle \;,\cr\cr
e_l |p\rangle & = \sum_{j=1}^l
c_{jl}^{-1} {P''_1(j,l;p_{jl}+1) P''_2(j,l;p_{jl}+1)
                               \over P''_3(j,l;p_{jl}+1) }
|p_{jl}+1 \rangle \;,\cr
}}
where $|p_{jl}\pm 1 \rangle $ denotes the state differing from
$|p\rangle $ by only $p_{jl}\rightarrow p_{jl}\pm 1$, and
\eqn\ePi{
\eqalign{
& P'_1(j,l;p)  = \prod_{i=1}^{l+1}
              [\varepsilon_{ij}(p_{i,l+1}-p_{j,l}+1)]^{1-\eta_{ijl}}
,\cr
& P''_1(j,l;p_{jl}+1)  = \prod_{i=1}^{l+1}
              [\varepsilon_{ij}(p_{i,l+1}-p_{j,l})]^{\eta_{ijl}}
,\cr}}
\eqn\ePii{
\eqalign{
& P'_2(j,l;p)  = \prod_{i=1}^{l-1}
              [\varepsilon_{ji}(p_{j,l}-p_{i,l-1})]^{\eta_{j,i,l-1}}
,\cr
& P''_2(j,l;p_{jl}+1)  = \prod_{i=1}^{l-1}
              [\varepsilon_{ji}(p_{j,l}-p_{i,l-1}+1)]^{1-\eta_{j,i,l-1}}
,\cr
}}
\eqn\ePiii{
\eqalign{
& P'_3(j,l;p)  = \prod_{{i=1 \atop i\ne j}}^{l}
            [\varepsilon_{ij}(p_{i,l}-p_{j,l})]^{1/2}
            [\varepsilon_{ij}(p_{i,l}-p_{j,l}+1)]^{1/2}
,\cr
& P''_3(j,l;p_{jl}+1)  = \prod_{{i=1 \atop i\ne j}}^{l}
            [\varepsilon_{ij}(p_{i,l}-p_{j,l}-1)]^{1/2}
            [\varepsilon_{ij}(p_{i,l}-p_{j,l})]^{1/2}
,\cr
}}
$\varepsilon_{ij}$ being the sign defined by
\eqn\eepsilon{\varepsilon_{ij}=\cases{ 1 & if $i\le j$\cr
                                      -1 & if $i>j$ . \cr}}

The parameters $\eta_{ijl}$ are introduced to break the symmetry
between the actions of $e_l$ and $f_l$, and to allow one to vanish
whereas the other does not. They will be taken to be $0$,
$1/2$ (the standard value) or $1$. They are not counted as
``continuous parameters'' in the following. These actions of the
generators on the Gelfand--Zetlin vectors define a module over $\uqn$
since they can formally be obtained from those of \rACi\ by a change
of normalization (See Appendix for details).

{}From the expression of the above matrix elements, it is obvious that
the indices $p_{il}$ belonging to the same line $l$ play a symmetric
role. They can formally be permuted. This remark did not hold in the
classical case when the indices were always related altogether by
triangular identities \etriangiden.

\medskip

Let us denote by $\alpha_i$ ($i=1,\cdots,N-1$) the simple roots of
$sl(N)$, and by $\alpha_{ij}=\alpha_i+\cdots +\alpha_{j-1}$ ($i<j$) the
positive roots. We define the raising generators
$e_{ij}\equiv e_{\alpha_{ij}}$ and $\tilde e_{ij}\equiv \tilde e_{\alpha_{ij}}$
for $i<j$ by
\eqn\eeij{
\cases{
e_{i,i+1} = \tilde e_{i,i+1}
\equiv e_i & for $i=1,...,N-1$ \cr
e_{i,j+1} = e_{ij}e_{j}-q^{-1}e_{j}e_{ij} & for $i<j$  \cr
\tilde e_{i,j+1} = \tilde e_{ij}e_{j}-q e_{j}\tilde e_{ij} & for $i<j$ ;  \cr
}}
The lowering generators $f_{ij}$ and $\tilde f_{ij}$ are defined by
the same inductions.

The action of these generators on the Gelfand--Zetlin representation
is given by

\eqn\ereprij{\eqalign{
f_{l,l+n+1} |p\rangle = &
\sum_{\scriptstyle j_i = 1,\cdots ,l+i \atop (i=0,\cdots, n)}
\left( \prod_{i=0}^n c_{j_i l+i}\right)
\left( \prod_{i=0}^{n-1} - \epsilon_{j_{i+1} j_i}\right)
q^{ - \left(p_{j_n l+n} - p_{j_0 l} + n \right) }      \cr
& \qquad \qquad
{\PP '_1(p) \PP '_2(p) \over \PP '_3(p) }
|p_{j_0 l}-1,\cdots, p_{j_n l+n}-1 \rangle \;,\cr\cr
e_{l,l+n+1} |p\rangle = &
\sum_{\scriptstyle j_i = 1,\cdots ,l+i \atop (i=0,\cdots, n)}
\left( \prod_{i=0}^n c_{j_i l+i} ^{-1} \right)
\left( \prod_{i=0}^{n-1} \epsilon_{j_{i+1} j_i}\right)
q^{ \left(p_{j_n l+n} - p_{j_0 l} \right) }   \cr
& \qquad \qquad
{\PP''_1(p_{j_i l+i}+1) \PP''_2(p_{j_i l+i}+1)
                               \over \PP''_3(p_{j_i l+i}+1) }
|p_{j_0 l}+1,\cdots, p_{j_n l+n}+1 \rangle \;,\cr
}}

For $\tilde f_{ij}$ and $\tilde e_{ij}$, just change the sign of the
exponent of $q$.

The symbols $\PP'_a$ and $\PP''_a$ ($a=1,\cdots,3$) denote the product
of all the factors coming from the product of $P'_a$ and $P''_a$
except those involving two of the modified indices $p_{j_i l+i}$
($i=0,\cdots,n$).
Hence, $\PP'_a$ (resp. $\PP''_a$) are the
common factors of the products
$f_l\cdots f_{l+n}$ (resp. $e_l\cdots e_{l+n}$) that
arise in the expansion of $f_{l,l+n+1}$ (resp. $e_{l,l+n+1}$), and
they do not depend on the order of the product. The
$q$-numbers involving differences of the indices $p_{j_i l+i}$
($i=0,\cdots,n$) (which depend on the order of the product),
gather and reduce to the single power of $q$.

\medskip
Note that for real \fp's of the indices,
each $q$-brackets in the above $P$ factors
is real. Thus, for example, if all the eponents $\eta_{ijl}$ are equal
to $1/2$, one has
hermitian conjugacy $e_\alpha = f_\alpha^\dagger$.
Moreover, $\left(k_i^{\pm 1}\right)^\dagger = k_i^{\mp 1}$.

\subsec{Restrictions on the values of the indices}

The above expressions \erepr\ involve
$P'_3$ and $P''_3$  in the denominators. They are well-defined if either
these denominators never vanish, or they are compensated by zeroes in
the numerators.

The first case is what happens for the most generic
representations
(i.e. generic indices $p_{il}$, not generic $q$),
with maximal dimension and number of parameters
\rDK. In this case, any
two indices of the same line $p_{il}$ and $p_{jl}$
have unequal
fractional parts, i.e. $\fp(p_{il} - p_{jl}) \neq 0$. Even after
translations by integers due to the action of the generators, the
$q$-numbers $[p_{il} - p_{jl}]$ never vanish and neither does the
whole
denominator.
The number of parameters and the dimension will be given in the
examples.

\medskip
On the other hand, some indices $p_{il}$ and $p_{jl}$
of the same line can have the same
fractional part, on the condition that zeroes in the numerators
compensate the denominators when they vanish. This happens if some
indices belonging to the adjacent lines $l\pm 1$
have the same fractional part as $p_{il}$ and $p_{jl}$.

We study here a sufficient condition for the matrix elements to be
well-defined. It is the condition that forbids to any pair of indices
of the same line to become equal under the action of the raising and
lowering generators.
This condition leads in particular to the case of usual
$q$-deformed representations.
Let us denote by $n_l(x)$ the number of indices $p_{il}$ of line $l$
with fractional part $x$.
One part of the condition is that this function
(which is non zero for only a
finite set of points of $\CC/{1\over 2}\ZZ$,
of course) obeys the following inequalities
\eqn\econdition{n_{l+1}(x) - 2 n_{l}(x) + n_{l-1}(x)\ge 0
\qquad \hbox{if} \qquad n_{l}(x)>1 \;.  }
Consider indeed two indices $p_{il}$ and $p_{jl}$  with the same
fractional part. The action of $f_l$ and $e_l$,
which is to translate them by
$\pm 1$, may make them become equal. If they reach the point where
$q^{\left(p_{il}\right)}=q^{\left(p_{jl}-1\right)}$,
then the denominator $P'_3$ or $P''_3$  vanishes in
the matrix element
$\langle p_{jl}-1|f_l|p\rangle$ or $\langle p_{il}+1|e_l|p\rangle$,
respectively. In order to keep these matrix elements finite,
some factors of
the numerator have to vanish also.
Furthermore, the matrix elements
$\langle p|e_l f_l|p\rangle$ and $\langle p|f_l e_l|p\rangle$
have to
remain finite, and moreover
$\langle p|e_l| p_{jl}-1 \rangle\langle p_{jl}-1|f_l|p\rangle$
and
$\langle p|f_l| p_{il}+1 \rangle\langle p_{il}+1|e_l|p\rangle$
have to be zero in order to
(i) keep the structure of module (preserve
    $[e_l,f_l]= {k_l -k_l^{-1} \over q-q^{-1}}$)
(ii) forbid  $q^{\left(p_{il}\right)}$ and
$q^{\left(p_{jl}\right)}$ to become equal, which would lead
to further divergences.
These constraints are satisfied if
\eqn\econdeta{\eqalign{
&\sum_{\left\{i' | 1\le i'\le l+1 \ \& \  q^{\left(p_{i',l+1}\right)}
=q^{\left(p_{jl}-1\right)} \right\}} (1-\eta_{i'jl})
+\sum_{\left\{i' | 1\le i'\le l-1 \ \& \  q^{\left(p_{i',l-1}\right)}
=q^{\left(p_{jl}\right)} \right\}} \eta_{j,i',l-1}
>{1\over 2}\cr
&\sum_{\left\{i' | 1\le i'\le l+1 \ \& \  q^{\left(p_{i',l+1}\right)}
=q^{\left(p_{jl}-1\right)} \right\}} \eta_{i'jl}
+\sum_{\left\{i' | 1\le i'\le l-1 \ \& \  q^{\left(p_{i',l-1}\right)}
=q^{\left(p_{jl}\right)} \right\}} (1-\eta_{j,i',l-1})
>{1\over 2}\cr
}}
plus analoguous constraints on $\eta_{i'il}$ and $\eta_{i,i',l-1}$.
The sum of conditions \econdeta\ implies that $q^{p_{il}}$ and $q^{p_{jl}}$
are separated, on
the discrete circle $\{q^{n p_{il}} \}_{n=0,...,m-1}$ by at least
two pairs of
indices. If the indices were not defined modulo $m$, we would need
at least one pair
$p_{i',l'}$, $p_{i'',l''}$, with $l',l''=l\pm 1$, with the same
fractional part as $p_{il}$ and $p_{jl}$, satisfying the usual
triangular identities
\eqn\eonepair{\cases{
p_{jl} >   p_{i' ,l' } \ge p_{il} & if $l' =l+1$  \cr
p_{jl} \ge p_{i' ,l' } >   p_{il} & if $l' =l-1$  \cr\cr
p_{jl} >   p_{i'',l''} \ge p_{il} & if $l''=l+1$  \cr
p_{jl} \ge p_{i'',l''} >   p_{il} & if $l''=l-1$  \cr
}}
Since the indices are periodic, $q^{p_{il}}$ and $q^{p_{jl}}$ have two
ways to
reach each other on the circle, so two pairs of indices belonging to
the adjacent lines are needed to prevent them from becoming equal, one
pair for each interval
$q^{p_{il}}$ and $q^{p_{jl}}$ define on the circle
$\{q^{n p_{il}} \}_{n=0,...,m-1}$.
\bigskip
\centerline{\epsfbox{uqsln2.eps}}
\medskip
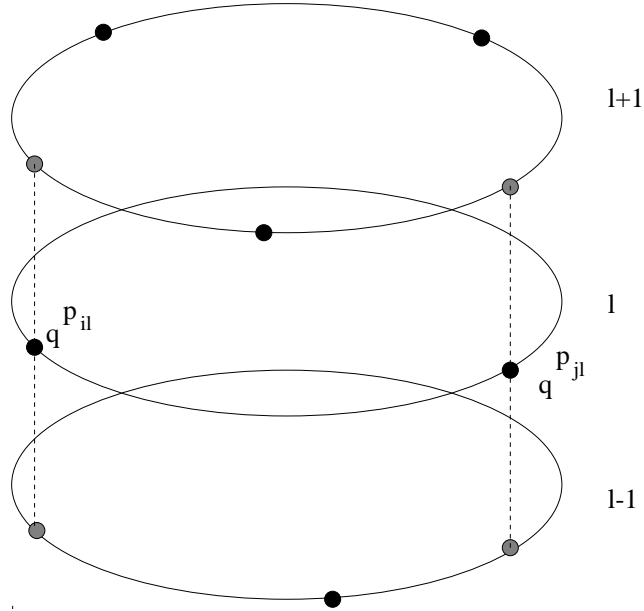
\centerline{\it Fig. 1: $q^{p_{il}}$ and $q^{p_{jl}}$ are separated, on each
side they define on the circle, }
\centerline{\it by two pairs of indices with same \fp\ from adjacent
lines $l-1$ or $l+1$.}
\bigskip
Then $n_{l}(x)$ indices on line $l$
with the same fractional part $x$ have to be separated by (at least)
$2 n_{l}(x)$ indices from lines $l\pm 1$ with the same fractional part
$x$, which is precisely the condition expressed by \econdition. If
$n_{l}(x)<2$, the above discussion of course does not apply and there
is no constraint at level $l$. If all the exponents $\eta$ are equal
to $1/2$, the condition \econdition\ is enough since \econdeta\ are then
automatically satisfied.

The set of indices with fractional part
$x$ can then be gathered into a sum of sub-triangles with a possible
line (which can be broken) starting from the lowest point of the
biggest triangle.
For this we use the symmetry among the indices that allows to reorder
each line. The triangular shape is the most natural since it recalls
the classical one with the triangular identities \etriangiden.
With respect to \rACi, the new point is the
possibility of a sum of triangles with the {\it same} fractional part,
and the line continuing the biggest triangle.
\bigskip
\centerline{\epsfbox{uqsln1.eps}}
\medskip
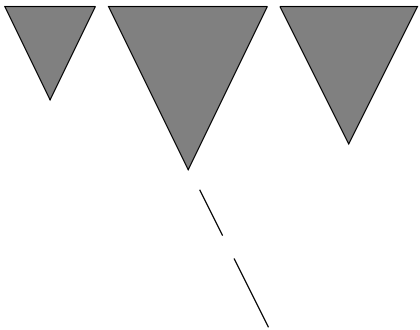
\centerline{\it Fig. 2: An admissible set of indices with the same \fp}
\bigskip

This analysis is done independently for each $x$, so that several sets
of indices with \fp\ $x$ that correspond to different $x$'s can
coexist.

A special case of sub-triangle of indices with same fractional part
has to be considered: when the indices $p_{jl}$ with
$1\le j \le l -N + N_1$ and
$N - N_1 + 1 \le l \le N$
(since the indices of the same line play the same role, we
choose here the upper left triangle for convenience)
satisfy the equalities
\eqn\etriangequal{p_{i,l+1} =  p_{il} = p_{i+1,l+1} + 1 \;.}
(similar to \etriangiden, but with equalities instead of
inequalities), then all these indices are frozen ($e_l$ or $f_l$
cannot  move them). All of them but $p_{1,N-N_1+1}=p_{1N}$ actually
disappear
from the matrix elements, since all numerators and denominators
involving them systematically cancel altogether (the exponents $\eta$
related to pairs of indices of the sub-triangle are chosen to be $1/2$).
Only the terms
$P'_1(j,N_1+1;p)$ and $P''_1(j,N_1+1;p)$ of the numerators still
involve $p_{1,N-N_1+1}=p_{1N}$.
In the actions of $e_l$ or $f_l$ \erepr,
the terms corresponding to shifts of
indices of this triangle can be removed.
This simplification allows us to forget completely about
the existence of
these indices, except $p_{1N}$ (fixed) of the first line, that
appears in $P'_1(j,N_1+1;p)$ and $P''_1(j,N_1+1;p)$.
The indices of the triangle are no more taken in account in the
function $n_l(x)$.
The existence of  $p_{1N}$ in $P'_1(j,N_1+1;p)$ and $P''_1(j,N_1+1;p)$
however changes the equation \econdition\ at level $N_1+1$ to
\eqn\econdfrozen{n_{N_1+2}(p_{1N}) - 2 n_{N_1+1}(p_{1N}) +
n_{N_1}(p_{1N}) + 1 \ge 0
\qquad \hbox{if} \qquad n_{l}(x)>1 \;.  }
By an abuse of notation, we write $n_l(p_{1N})$ instead of
$n_l(\fp(p_{1N}))$. This should not lead to any confusion.

Again, several such sub-triangles can coexist.

\newsec{Examples}

\subsec{Periodic (cyclic) representations}
The most generic representations (do not confuse with generic $q$)
are, as explained above, those for which $n_l(x)<2$ for all $l$ and
$x$, i.e. two indices of the same line do not have the same fractional
part.
The indices are also not bounded from above or below by other indices
of adjacent lines. The dimension is the maximal allowed dimension
when $q$ is a $m$-th root of unity, i.e., for $m$ odd,
$(m)^{N(N-1)/2}$.
Each of the $N(N-1)/2$ indices $p_{il}$ ($l<N$) take $m$ values.
(For $m$  even, case we do not consider here,
the representation would not be irreducible unless we
identify $|p_{il}+m/2,p_{jl}+m/2\rangle$. In this case,
the dimension is given in
\rACi.)

These representations are called periodic (or cyclic) since for
$\alpha$ a positive root,
$f_\alpha^{m}$ and  $e_\alpha^{m}$ act as (generally non zero)
scalars on them.

The continuous parameters for periodic representations are
\item{-- } the $c_{jl}$'s, for $1\le j\le l <N$
\item{-- } the $p_{jl}$'s, for $1\le j\le l \le N$ (in fact only by
their $q$-th power, and modulo integer powers of $q$ for those which
do not belong to the top line).

The total number of parameters is then $N^2-1$ (after taking in account
the constraint $\sum p_{iN}=0$ or $p_{NN}=0$).
These parameters are indirectly related to the values of the $N^2-1$
central operators $f_\alpha^{m}$, $e_\alpha^{m}$ and $k_i^{m}$. The
values of the $q$-deformed ordinary Casimir operators are actually not
independent of those (see Sect. 4).
Both the dimension and
number of parameters agree with \rDK.
The values of $\eta_{ijl}$ do not matter in this case.

\medskip
All the other examples correspond to less generic cases.
For non-generic representations, the parameters live on
sub-manifolds of the $N^2-1$
dimensional manifold of the whole set of parameters (for instance
fractional parts of some indices become equal). This possibly leads to
\item{-- } loss of periodicity of some $f_\alpha$ or  $e_\alpha$.
\item{-- } reduction of dimension.

\subsec{Semi-periodic (semi-cyclic) representations}

Semi-periodic representations are highest weight representations for
which the lowering generators are periodic (resp. lowest weight
representations with periodic raising generators).

Take $n_l(x)\le 1 \quad \forall l,x$, and $\fp(p_{il} - p_{i,l+1})=0$ for
$i\le l$. Choose $\eta_{ijl}=1 \quad \forall i,j,l$.

Then we get a highest weight representation. The highest weight state
$|p_0\rangle$ is given by $p_{il}=p_{iN} \quad \forall i\le l\le N$.
On this state,
\eqn\ehwv{e_l |p_0\rangle=0 \qquad \forall l<N \;,
\qquad \hbox{and hence} \qquad e_\alpha |p_0\rangle=0
}
for all the raising generators $e_\alpha$.

The dimension of these representations is the same as for periodic
representations, and the number of parameters is $(N-1)(N+2)/2$
(i.e. $N-1$ independent $p_{iN}$ and $N(N-1)/2$ $c_{jl}$'s). The
$f_\alpha$'s remain periodic on these representations;
for this reason we call them semi-periodic. The
values of the central operators $f_\alpha^m$ and $k_i^m$
are independent and
related to the remaining parameters. Note that the vanishing of a
$f_\alpha^m$ is not directly related to particular equality of some
\fp\ of indices, but rather on more general algebraic equations among the
$q^{m p_{il}}$ and the $c_{jl}^m$.

Semi-periodic representations can also be lowest weight
representations, if $\eta_{ijl}=0$. More complicated examples with
mixed vanishing of $e_\alpha^m$ and $f_\beta^m$ exist, which can be
obtained directly by a suitable choice of the parameters, or also by
braiding action of the Weyl group \rLusJAMS\
on a highest weight semi-periodic representation.

This example of representations could not be taken in account our
first approach \rACi, because the symmetry between raising and
lowering generators was not enough broken.

\subsec{Nilpotent representations}

Nilpotent representations are representations with a highest weight
vector and a lowest weight vector, and hence nilpotent action of all
the raising and lowering generators. They still have complex
parameters related to the values of the operators $k_l^{m}$.

Take as before
$n_l(x)\le 1 \quad \forall l,x$, and $\fp(p_{il} - p_{i,l+1})=0$ for
$i\le l$. Choose now $\eta_{ijl}={1\over 2} \quad \forall i,j,l$.

The dimension is in this case
$(m)^{N(N-1)/2}$, i.e. the same as for periodic representation.
This representation is nilpotent (i.e. $f_\alpha^{m}=e_\alpha^{m}=0$
for every positive root $\alpha$) and it is not necessary to consider
the indices $p_{il}$ modulo $m$, since the range of values for each
index is bounded above and below by indices of adjacent lines, the
upper ones being fixed.
The nilpotent representation is characterized by $N-1$
parameters (the $p_{iN}$'s with $\sum p_{iN}=0$ or $p_{NN}=0$),
corresponding the values of the operators $k_l^{m}$, or
to the $q$-deformed usual Casimir operators. (The parameters
$c_{jl}$ can here be set to one by a change of normalization.)

We have the same highest weight vector as for semi-periodic
representations. But the   $f_\alpha$'s are no longer periodic. The
specification
$p_{il}=p_{i,l+1}-m+1 \quad \forall i\le l\le N$
indeed defines the lowest weight vector of this representation.

In \rACi, we did not have the correct number of parameters for this
kind of representations.

\subsec{Usual $q$-deformed representations}

We now consider usual representations, i.e. those that are the
$q$-deformations of the classical representations, in the limit where
$q$ is a root of unity. In the classical case, or when $q$ is not a
root of unity, the Gelfand--Zetlin indices are integers and there is no
reason to define them modulo $m$. When we take the limit where
$q$ is a root of unity, we expect no periodicity of the indices
(usual representations are highest-weight and lowest-weight ones) so
they are not considered modulo $m$.

The usual representations correspond in the Gelfand--Zetlin formalism
\rJimGZ\ to usual choice where $n_l(0) = l$ for
all $l\le N$, i.e. all the indices $p_{il}$ have the same fractional part
$0$. The exponents $\eta_{ijl}$ are set to $1/2$.
{\sl No continuous parameter} survives since the $c_{jl}$'s can be
absorbed in a change of normalization.

Only a finite number of representations, those with a highest weight
satisfying \refs{\rDob, \rFGP, \rACi}
\eqn\elimit{p_{1N} - p_{NN} \le m \;.}
are well-defined in the Gelfand--Zetlin formalism for $q^m=1$.
The condition \elimit\ expresses the fact that the $q$-th powers of the
indices of the first line do not wind more than exactly once
around the circle $\{q^n \}_{n=0,...,m-1}$.
This condition is also the unitarity condition, i.e.
all the matrix elements of $e_l$ and $f_l$ are real for the usual
representations, and the matrices of $k_l$ are unitary. Furthermore,
the matrices representing $e_l$ is the transposed matrix of that
of $f_l$.

The representations with highest weights that do not obey \elimit\
are considered in the following subsection,
although they are not all atypical.

\subsec{Atypical representations}

We consider here the quantum analogue of classical (highest weight and
lowest weight) irreducible
representations with a highest weight that does not obey \elimit.
When $q^m=1$ these representations are not always irreducible, since
some new singular vectors arise in the corresponding Verma modules
\rDob, that are not obtained from the highest weight vector by
action of the translated Weyl group.
Quotienting by the sub-representation generated by these
singular vector leads to new irreducible representations that we call
atypical by analogy with the case of superalgebras.

The Gelfand--Zetlin basis in the form we consider is not yet totally
adapted for atypical representations. This has to be compared with the
fact that, for superalgebras, the atypical
representations are more difficult to describe with the
Gelfand--Zetlin than the typical ones:
the atypical
representations of some superalgebras or quantum superalgebras
were obtained, for example, in  \refs{\rPal,\rPalTol} in the case of
$gl(n|1)$ and in \rNAKyStoi\ in the case of $\cU_q(gl(2|2))$,
but the general case is not yet written.

It seems here that a further adaptation of the Gelfand--Zetlin basis
to the atypical case is possible. We already obtained some
examples of (reducible or irreducible)
representations that do not obey \elimit. A general study of this case
will be the subject of another work. Note that the formalism of \rBKW,
in which the matrix elements do not contain divergences, provides the
atypical representations of $\cU_q(sl(3))$.

Some atypical representation can also be obtained as degenerations of
periodic representations, by taking the appropriate limit of the
parameters.
Consider as explained before the possibility of ``freezing'' a
sub-triangle of indices
$p_{jl}$ with
$1\le j \le l -N + N_1$ and
$N - N_1 + 1 \le l \le N$.
Remember that these indices are not taken in account in
$n_l(p_{1N})$, and that the inequality \econdfrozen\ holds instead of
\econdition\ for $x=\fp(p_{1N})$ and $l=N_1+1$.

Choosing $n_l(p_{1N})=0$ for $l=1,...,N_1$ leads only to representations
described in the next subsection as ``partially periodic''. Some of
the generators
remain indeed periodic. With $N_1=N-1$, they are
``flat'', i.e. the multiplicities of their weights are always $1$
\refs{\rACi,\rACii}.

Choosing $n_l(p_{1N})\neq 0$ can however lead to some cases of
atypical representations.
In particular, with $n_l(p_{1N}) = 1$ and $N_1=N-1$,
we get truncated flat representations \rAsutroisq\ that can also be
seen as an irreducible part of the limit when $q^{m}=1$ of
representations with $p_{1N}-p_{NN}=m+1$, i.e. just after the limit
given by \elimit\ \refs{\rDob, \rBKW, \rFGP}.

We can recall as examples the cases of the representations of
$\cU_q(sl(3))$ of dimensions $7$ for $m=3$ and $18$ or $19$ for
$m=5$.
A formula for the dimensions of these representations in the
case of $\cU_q(sl(3))$ is, with $N_1=2$ and $p_{33}=0$,
\eqn\edim{\eqalign{
d & ={m}^2 - d_1 - d_2 \qquad \cases{
            d_1 = {1\over 2} p_{13}(p_{13}-1) \cr
            d_2 = {1\over 2} (m-p_{13}+1)(m-p_{13}) \cr} \cr
  & = d'_1 -d'_2 \qquad \cases{
            d'_1 = {1\over 2} (m+1)p'_{23}(m+1-p'_{23}) \cr
            d'_2 = {1\over 2} (m-1)(p'_{23}-1)(m-p'_{23}) \cr} \cr
}}
The first expression corresponds to the truncation of the flat
representation of dimension ${m}^2$ by its
two triangular sub-factors of dimensions $d_1$ and $d_2$
(See the figure in \rAsutroisq). The
second expression corresponds to the
same representation seen as the irreducible part of the limit when
$q^{m}=1$  of the representation of dimension $d'_1$
with first line of indices
(or highest weight) $|p'_{13}=m+1,p'_{23}=p_{13},p'_{33}=0\rangle$
violating
\elimit\ by $1$; $d'_2$ being the dimension of its sub-representation
characterized by
$|p'_{13}=m,p'_{23}=p_{13},p'_{33}=1\rangle$. The second expression is
a particular case of those classified in \rDob.

The generalization of these cases to $\uqn$ with $N_1=N-1$ is
straightforward (flat representations).
Different values for $N_1$ provide
other interesting examples.

\subsec{Partially periodic representations}

First note that, since the whole set of indices $p_{il}$ is
defined up to an
overall constant, the case of usual $q$-deformed representations
can be written with $n_l(x) = l$
for any given value $x\in \CC$ instead of $x=0$.

As in \rACi\ (see Sect. 3),
one can put in some sub-triangles (those defined by
\econdition\ and Fig. 2) of sizes $N_1$, ..., $N_a$,
with $N_1+...+N_a\le N$,
the indices
corresponding to usual representations
of some $\cU_q(sl(N_1))$, ... ,
$\cU_q(sl(N_a))$.

This prescription reduces the dimension with respect to the maximal
one. Each triangle indeed contributes to the dimension by a factor
equal to the dimension of the related usual representation of the
corresponding $\cU_q(sl(N_i))$,
instead of a factor equal to $(m)^{N_i(N_i-1)/2}$.
The number of parameters is also reduced, since all the indices
of a given sub-triangle have the same fractional part, whereas the
corresponding $c_{jl}$'s are $1$.

\medskip

The atypical representations of smaller
$\cU_q(sl(N_i))$
can also be used as the usual ones to
construct partially periodic representations $\uqn$.

\newsec{Application: a set of relations in the centre of $\uqn$}

The centre of $\uqn$ is generated by the operators
$f_\alpha^{m}$, $e_\alpha^{m}$, $k_i^{m}$, and the $q$-deformed
classical Casimirs $\cC_i$ \rDK.
Let us introduce supplementary Cartan generators $k_{\pm \epsilon_i}$
for $i=1,...,N$,
that are needed to
write the $q$-deformed Casimirs.
These generators are such that
$ k_{\epsilon_i} k_{-\epsilon_{i+1}} = k_{\alpha_i} \equiv k_i $ and
$\prod_{i=1}^{N} k_{\epsilon_i} =1$.
They satisfy the relations
\eqn\erelkepsilon{
\eqalign{
  k_{\epsilon_i}  e_{j}  k_{\epsilon_i}^{-1}
&= q^{\delta_{ij} -\delta_{i-1, j}} e_{j} \;,
\cr
  k_{\epsilon_i}  f_{j}  k_{\epsilon_i}^{-1}
&= q^{-\delta_{ij} +\delta_{i-1, j}} f_{j}  \;  .
}}
By convention, $k_{\zeta+\xi}=k_\zeta k_\xi$.
The operators  $k_{\pm \epsilon_i}^{m}$ are also central.

A set of generators for the deformed classical
centre (i.e. the whole centre
when $q$ is not a root of unity) is given by
\eqn\eCiHC{\left\{\cC _i=h^{-1} \left(
\sum_{1 \leq j_1 < \cdots < j_i \leq N}
k_{2\epsilon_{j_1}} \cdots
k_{2\epsilon_{j_i}}
\right)\right\}_ {i=1,...,N-1}\;\;,}
where $h$ is the Harish--Chandra isomorphism \refs{\rRosHC, \rDKP}\
between the $q$-deformed classical
centre and the algebra of symmetric
polynomials in the $k_{2\epsilon_i}$'s.

This isomorphism $h$ can be written as $h=\gamma^{-1} \circ h'$,
with the following notations:
$h'$ is the projection on $\cU^0$, within the direct sum
$\cU=\cU^0 \oplus (\cU^-\cU + \cU \cU^+)$, with $\cU\equiv\uqn$, and
where $\cU^0$ (resp. $\cU^+$ and
$\cU^-$) is the sub-algebra of $\uqn$ generated
by the $k_{\pm \epsilon_i}$'s (resp. $e_i$'s and $f_i$'s).
$\gamma$ is the
automorphism of $\cU^0$ given by
$\gamma(k_{2\epsilon_i}) = q^{N+1-2i} k_{2\epsilon_i} $.

Let us write $\cC_i$ as a function (actually a non commuting
polynomial) of the parameter $q$ and the generators
\eqn\eCifun{ \cC _i= \cF_i \left( q,k_{2\epsilon_j},
         \lambda e_\alpha, \lambda f_\alpha \right)
\qquad
   j=1,...,N, \qquad \alpha\in{\hbox{set of positive roots}}
}
with $\lambda = q-q^{-1}$.

Then the comparison of the actions of the $m$-th powers of the
generators with the actions of the generators themselves
on periodic representation, provides us the relations that hold in the
centre of the algebra:
\eqn\erelations{\cP_{i,m}^{(N)}  (\cC _1,...,\cC _{N-1}) =
\cF_i \left( q^{m}=1 , k_{2\epsilon_j}^{m} ,
         \lambda^{m} e_\alpha^{m}, \lambda^{m} f_\alpha^{m}
\right)\;,
}
where $\cP_{i,m}^{(N)}$ is the polynomial such that
\eqn\ePolynome{\cP_{i,m}^{(N)}  (h(\cC _1),...,h(\cC _{N-1})) =
 \sum_{1 \leq j_1 < \cdots < j_i \leq N}
k_{2\epsilon_{j_1}}^{m} \cdots
k_{2\epsilon_{j_i}}^{m}\;.}
(See \rAB\  for details and a proof of these relations for $i=1$ or
$N-1$.
A proof of these relations for $i=2, \cdots N-2$ will be given
elsewhere.)

In \erelations, the left hand side is a polynomial in the $q$-deformed
classical Casimirs, whereas the right hand side is a function of
operators that are central  only when $q$ is a root of unity. The nice
feature is that this
function is, up to numerical coefficients, the same as the polynomial
that defines the $i$-th Casimirs in terms of the generators.

\newsec{Appendix}

We give here the transformation that relates the representation given
by \erepr--\ePiii\ to the standard $q$-deformed Gelfand--Zetlin basis
\rJimGZ.
Let us denote by $|p)$ be the vectors of the representation of $\uqn$
in the Gelfand--Zetlin basis with all the exponents $\eta$ equal to
$1/2$, as defined in \rJimGZ. For generic values of the indices
$p_{il}$, we define the new basis
\eqn\etransformation{
\eqalign{
& |p> = \lambda(p) |p) \cr
& \lambda(p) = \prod_{l=1}^{N-1} \prod_{j=1}^l \lambda_{jl}(p) \cr
& \lambda_{jl}(p) = \prod_{i=1}^{l+1}
       \left( \Gamma_q \left( \epsilon_{ij}(p_{i,l+1}-p_{jl}+1/2)+1/2
        \right)\right)^{\epsilon_{ij}(\eta_{ijl}-{1\over 2})} \cr
&\phantom{\lambda_{jl}(p) = }
     \prod_{i=1}^{l-1}
     \left( \Gamma_q \left( \epsilon_{ji}(p_{j,l}-p_{i,l-1}+1/2)+1/2
        \right)\right)^{\epsilon_{ji}(\eta_{j,i,l-1}-{1\over 2})}\;,
}}
where the function $\Gamma_q$
obeys
\eqn\egamma{\Gamma_q(x+1) = [x] \Gamma_q(x)\;.}
This function is a straightforward
adapation of the definition of \rGasperRahman, in which the
definition of $q$-number is different.
This transformation is well-defined when the \fp's of the indices
$p_{il}$ are unequal. It works formally on infinite dimensional
representations with no identification of the indices modulo $m$.
The actions of the generators on this basis are given by \erepr, which
then defines a module on $\uqn$. Quotienting then by the
identification of the indices modulo $m$, and exploration of the
parameter manifold leads to the examples of
representations described in this paper.

\bigskip
{\bf Acknowledgements:} We would like to thank M. Bauer for numerous
discussions on the relations in the centre.
This work is supported in part by the EEC contracts No. SC1-CT92-0792
and No. CHRX-CT93-0340.

\vfill\eject
\baselineskip=16pt
\listrefs

\bye